\documentclass[12pt]{article}
\usepackage{tcolorbox}
\usepackage{amssymb,amsmath,bm}
\usepackage[colorlinks=true, urlcolor=blue,linkcolor=blue, citecolor=blue]{hyperref}
\usepackage{mathrsfs,verbatim}
\usepackage{caption,cleveref }
\usepackage{graphicx, float}
\usepackage{constants,color,appendix,  mathtools, tikz, xcolor}
\usepackage{enumerate}
\usepackage{bbm}
\usepackage{tkz-euclide}
\usepackage{tikz}
\usetikzlibrary{trees}
\usepackage{pdfpages}
\usepackage{xcolor}

\usetikzlibrary{calc}
\usepackage{lipsum}

\newcommand{\Pro}{\mathbb{P}}

\usepackage{appendix}
\newcommand{\ds}{\displaystyle}

\begin{document}
\newcommand{\bea}{\begin{eqnarray}}
\newcommand{\ena}{\end{eqnarray}}
\newcommand{\beas}{\begin{eqnarray*}}
\newcommand{\enas}{\end{eqnarray*}}
\newcommand{\beq}{\begin{equation}}
\newcommand{\enq}{\end{equation}}
\def\qed{\hfill \mbox{\rule{0.5em}{0.5em}}}
\newcommand{\bbox}{\hfill $\Box$}
\newcommand{\ignore}[1]{}
\newcommand{\ignorex}[1]{#1}
\newcommand{\wtilde}[1]{\widetilde{#1}}
\newcommand{\mq}[1]{\mbox{#1}\quad}
\newcommand{\bs}[1]{\boldsymbol{#1}}
\newcommand{\qmq}[1]{\quad\mbox{#1}\quad}
\newcommand{\qm}[1]{\quad\mbox{#1}}
\newcommand{\nn}{\nonumber}
\newcommand{\Bvert}{\left\vert\vphantom{\frac{1}{1}}\right.}
\newcommand{\To}{\rightarrow}
\newcommand{\supp}{\mbox{supp}}
\newcommand{\law}{{\cal L}}
\newcommand{\Z}{\mathbb{Z}}
\newcommand{\mc}{\mathcal}
\newcommand{\mbf}{\mathbf}
\newcommand{\tbf}{\textbf}
\newcommand{\lp}{\left(}
\newcommand{\limm}{\lim\limits}
\newcommand{\limminf}{\liminf\limits}
\newcommand{\limmsup}{\limsup\limits}
\newcommand{\rp}{\right)}
\newcommand{\mbb}{\mathbb}
\newcommand{\rainf}{\rightarrow \infty}
\newtheorem{problem}{Problem}[section]
\newtheorem{exercise}{Exercise}[section]
\newtheorem{theorem}{Theorem}[section]
\newtheorem{corollary}{Corollary}[section]
\newtheorem{conjecture}{Conjecture}[section]
\newtheorem{proposition}{Proposition}[section]
\newtheorem{lemma}{Lemma}[section]
\newtheorem{definition}{Definition}[section]
\newtheorem{example}{Example}[section]
\newtheorem{remark}{Remark}[section]
\newtheorem{solution}{Solution}[section]
\newtheorem{case}{Case}[section]
\newtheorem{condition}{Condition}[section]
\newtheorem{defn}{Definition}[section]
\newtheorem{eg}{Example}[section]
\newtheorem{thm}{Theorem}[section]
\newtheorem{lem}{Lemma}[section]
\newtheorem{soln}{Solution}[section]
\newtheorem{propn}{Proposition}[section]
\newtheorem{ex}{Exercise}[section]
\newtheorem{conj}{Conjecture}[section]
\newtheorem{pb}{Problem}[section]
\newtheorem{cor}{Corollary}[section]
\newtheorem{rmk}{Remark}[section]
\newtheorem{note}{Note}[section]
\newtheorem{notes}{Notes}[section]
\newtheorem{readingex}{Reading exercise}[section]
\newcommand{\pf}{\noindent {\bf Proof:} }
\newcommand{\proof}{\noindent {\it Proof:} }
\frenchspacing

\tikzstyle{level 1}=[level distance=2.75cm, sibling distance=5.65cm]
\tikzstyle{level 2}=[level distance=3cm, sibling distance=2.75cm]
\tikzstyle{level 3}=[level distance=3.9cm, sibling distance=1.5cm]

\tikzstyle{bag} = [text width=10em, text centered] 
\tikzstyle{end} = [circle, minimum width=3pt,fill, inner sep=0pt]


\title{\bf A study on random permutation graphs}
\author{O\u{g}uz G\"{u}rerk\footnote{Bo\u{g}azi\c{c}i University, Department of Mathematics, Istanbul, Turkey, e-mail:  oguz.gurerk@boun.edu.tr} \hspace{0.2in} \"{U}m\.{i}t I\c{s}lak\footnote{Bo\u{g}azi\c{c}i University, Department of Mathematics, Istanbul, Turkey, e-mail: umit.islak1@boun.edu.tr}  \hspace{0.2in} Mehmet Akif Y{\i}ld{\i}z\footnote{University of Amsterdam, Korteweg-de Vries Institute for Mathematics, Amsterdam, the Netherlands, e-mail: m.a.yildiz@uva.nl}} \vspace{0.25in}

\maketitle

 
\begin{abstract}
 For a given permutation $\pi_n$ in $S_n$, a random permutation graph is formed by  including  an edge between two vertices $i$ and $j$ if and only if $(i - j) (\pi_n(i) - \pi_n (j)) < 0$. In this paper, we study various  statistics of random permutation graphs. In particular, the degree of a given node, the number of nodes with a given degree, the number of isolated vertices, and the number of cliques are analyzed. Further, explicit formulas for the probabilities of  having a given number of connected components and isolated vertices  are obtained.

\bigskip



\end{abstract}

\section{Introduction}\label{sec:intro}

Let $S_n$  be the set of all permutations on $[n] := \{1,2,\ldots,n\}$.  For a given permutation $\pi_n$ in $S_n$, we may form a corresponding graph by including  an edge between the given  two vertices $i$ and $j$ if and only if $(i - j) (\pi_n(i) - \pi_n (j)) < 0$. In other words, there turns out to be an edge between two vertices if and only if the corresponding values in the permutation form an inversion. The graph resulting from this process is known to be a \emph{permutation graph}. If we consider the uniform distribution over $S_n$, then the resulting model is known to be the \emph{random permutation graph model}. 
 The purpose of current manuscript is to study certain degree related statistics in this class of random permutations.

Permutation graphs were first introduced and discussed by Even, S., Lempel, A. and Pnueli, A.   in  \cite{ELP:1972} and  \cite{PLE:71}. These two cited works provide a characterization for being a permutation graph in terms of transitive orientability. In particular, it is shown that a graph is a permutation graph if and only if both the graph and its complement  are transitively orientable, which means that whenever the edges are assigned directions, the presence of the edges $(x,y)$ and $(y,z)$ implies the existence of the edge $(x,z)$. It is   known that testing whether a graph is a permutation graph can be done in linear time.  Let us also note that in these two papers, they also provided a polynomial time algorithm to find a transitive orientation when it exists. 
One particular interest in permutation graphs stems from the fact that they are perfect graphs. This enhances their applications in several fields such as channel routing, scheduling, memory allocation, genomics and bioinformatics.  See \cite{BP:96} as an exemplary work.

In \cite{AH:2016}, they emphasize that various features of random permutation graphs can be understood by turning the problem into a random permutation problem. That will be the path we   follow below.   Most of the statistics we will work on will be  related to the descent structure of the underlying permutation. We refer to \cite{Bona:2016} as  a general reference for combinatorics of permutations.  \cite{BDJ:1999}, \cite{DFG:2008} and \cite{O:2000} are some important works on statistics of random permutations, but the field is so huge that it is hard to include enough references here. 


 We will now  briefly sketch what we discuss below.  First,  let us note that from here on  we will be using one line notation for permutations; for example, $\binom{1 \hspace{0.1in}  2 \hspace{0.1in}  3 \hspace{0.1in}  4}{3 \hspace{0.1in}  1 \hspace{0.1in}  2 \hspace{0.1in}  4} = (3,1,2,4)$.     $\mc{G}_n$ will be denoting  a random permutation graph on $n$ vertices throughout the paper. The following figure shows a sample from this model when $n = 5$.

 \begin{figure}[H]
 \begin{center}
 \includegraphics[scale=0.4]{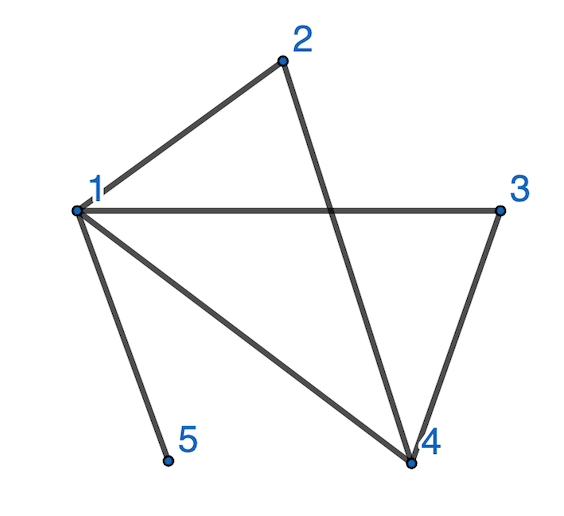}
 \caption{The permutation graph corresponding to the permutation $\binom{1 \hspace{0.1in}  2 \hspace{0.1in} 3 \hspace{0.1in} 4 \hspace{0.1in} 5}{5 \hspace{0.1in} 2 \hspace{0.1in} 3 \hspace{0.1in} 1 \hspace{0.1in} 4}$.}
 \end{center}
 \end{figure}
 
 Our first  focus in Sections \ref{sec:asymid} and \ref{sec:degothnod} will be on the asymptotic distribution of the degree of a given node. In particular,   we provide a very simple proof for a central limit theorem for the mid-node which was previously proved in \cite{BM:2017} by using different techniques. We do not restrict ourselves here   just to   the mid-node, but we also prove a central limit theorem for any given fixed node $k$ as the number of nodes  $n$ grows. 

Our interest afterwards is on the number of isolated vertices. This first requires an understanding of connected components which is done in Section \ref{sec:conncom}. The content of this section is somewhat already available in the literature, but we will be including all the details for the sake of completeness. Regarding isolated vertices, we study   the probability that the resulting graph has exactly some given number of isolated vertices (Section \ref{sec:numiso}), and also the expected  number of such vertices (Section \ref{sec:ndegd}).  Our treatment on the expectation of number of  isolated vertices will be  more general, and we will derive a formula for the expected number of degree $d$ vertices in a random permutation graph covering the number of isolated vertices, leaves, etc., as special cases. 

Later in Section \ref{sec:clique}, we discuss the number of $m$-cliques $K_m$, $m \in\mbb{N}$ fixed, and     prove that the following central limit theorem holds by making use of the theory of $U$-statistics: 
$$\frac{K_m - \mathbb{E}[K_m] }{\sqrt{Var(K_m)}}
\longrightarrow_d \mathcal{Z}, \qquad \text{as} \; \; n\rightarrow \infty,$$    where $\mc{Z}$ is a standard normal random variable.  
It is then easily  seen that the same central limit theorem also holds  for the number of cycles in $\mc{G}_n$ of  length at least $m$, discussed in the same section.  
 


We  fix some notation now. First, $=_d$ and $\rightarrow_d$ are used for equality in distribution and convergence in distribution, respectively. $\mathcal{Z}$ denotes the standard normal random variable. Finally, for two sequences $a_n$ and $b_n$, we write $a_n \sim b_n$  if and only if  $\lim_{n\rightarrow \infty} a_n / b_n =1$.


\section{Asymptotic distribution of the  degree  of mid-node}\label{sec:asymid}

Let $\mc{G}_n$ be  a random permutation graph on $n$ vertices. The purpose of this section is to prove that the mid-node  in $\mc{G}_n$  satisfies a central limit theorem.   In order to  avoid expressions involving  the floor function, we  assume  here that   $n$ is even.     Our  result  in this section   was also proven in \cite{BM:2017},  Theorem   3.4   in cited paper. However, we believe that our treatment is much   elementary and  has a more probabilistic flavour.  One other advantage is that our approach extends to other nodes as we shall see in  Section \ref{sec:degothnod}.

\begin{theorem}
Let $\mc{G}_n$ be  a random permutation graph on $n$ vertices. Then, the degree $d(n/2)$ of the mid-vertex satisfies 
\begin{equation}\label{eqn:cltwu}
\frac{d(n/2) - n / 2}{\sqrt{n}} \rightarrow_d \mc{N}(0, U(1-U))
\end{equation}
 where $U$ is a uniform random variable over $(0,1)$ independent of $\mc{G}_n$ and $\mc{N}(0,U(1-U))$ is a normal random variable with the corresponding (random) parameters. 
\end{theorem}

\pf
The expectation and variance will be special cases of the results presented in the next section (Theorem \ref{thm:expvardeggen}), so we skip these computations here. 

Now, observe that we may write $d(n/2)  =  \sum_{j=1}^{n/2 -1 } \mbf{1}(\pi_j > \pi_{n/2}) + \sum_{j = n/2 +1}^n \mathbf{1}(\pi_{n/2} > \pi_j)$. 
 Letting then   $X_1,X_2,\ldots$ be an i.i.d. sequence of uniform random variables over $(0,1)$, by a result often attributed to R\'enyi, we may express $d(n/2)$ as the following: 
 $$d(n/2) =_d \sum_{j=1}^{n/2 -1 } \mbf{1}(X_j > X_{n/2}) + \sum_{j = n/2 +1}^n \mathbf{1}(X_{n/2} > X_j).$$
 Let $X_1^*,X_2^*,\ldots$ be an independent copy of  $X_1,X_2,\ldots$, and let $\alpha\in \mathbb{R}$. We then have 
 
\begin{eqnarray*}
\mbb{P} \lp \frac{d(n/2) - n/2	}{ \sqrt{n}} \leq \alpha \rp &=& \mbb{P} \lp \frac{\sum_{j=1}^{n/2 -1 } \mbf{1}(X_j > X_{n/2}) + \sum_{j = n/2 +1}^n \mathbf{1}(X_{n/2} > X_j) - \frac{n}{2}}{\sqrt{n}} \leq \alpha \rp  \\
&=& \int_0^1 \mbb{P} \lp \frac{\sum_{j=1}^{n/2 -1 } \mbf{1}(X_j > u) + \sum_{j = n/2 +1}^n \mathbf{1}(u > X_j) - \frac{n}{2}}{ \sqrt{n}} \leq \alpha \rp  du \\ 
&=&   \int_0^1 \mbb{P} \lp \frac{\sum_{j=1}^{n/2 -1 } \mbf{1}(X_j > u) + \sum_{j = n/2 +1}^n(1 -  \mathbf{1}(u<  X_j)) - \frac{n}{2}}{ \sqrt{n}} \leq \alpha \rp  du \\ 
&=&   \int_0^1 \mbb{P} \lp  \frac{\sum_{j=1}^{n/2 -1 } \mbf{1}(X_j > u) - \sum_{j = n/2 +1}^n \mathbf{1}(u < X_j)}{ \sqrt{n}} \leq \alpha \rp  du \\ 
&=&   \int_0^1 \mbb{P} \lp  \frac{\sum_{j=1}^{n/2 -1} \mbf{1}(X_j > u) - \sum_{j = 1}^{n/2} \mathbf{1}(u < X_j^*)   }{ \sqrt{n} } \leq \alpha \rp  du \\ 
&=&  \int_0^1 \mbb{P} \lp  \frac{\lp \sum_{j=1}^{n/2 -1 } ( \mbf{1}(X_j > u) - \mathbf{1}(X_j^* > u)) \rp - \mathbf{1}(X_{n/2}^* > u) }{ \sqrt{n}} \leq \alpha \rp  du \\ 
\end{eqnarray*}
Now, the probability in the integrand for a given  $u$ converges to $\int_{- \infty}^{\alpha} \frac{1}{\sqrt{2 \pi} \sqrt{u ( 1 - u)}} e^{-\frac{x^2}{2 \sqrt{u (1 - u)}}} dx $ by using the standard  central limit theorem for i.i.d. random variables. Integration over $u$ from 0 to 1,  together with an application of dominated convergence theorem, yields as $n\rightarrow\infty$, $$\mbb{P} \lp \frac{d(n/2) - n/2}{ \sqrt{n}} \leq \alpha \rp  \rightarrow   \int_0^1 \int_{- \infty}^{\alpha} \frac{1}{\sqrt{2 \pi} \sqrt{u ( 1 - u)}} e^{-\frac{x^2}{2 \sqrt{u (1 - u)}}} dx du.$$ But the last expression equals $\mathbb{P}(N(0,U(1-U)) \leq \alpha)$.   The asserted convergence in distribution follows. 
\bbox

\section{Degrees of other  nodes}\label{sec:degothnod}

Previous section  focused on the mid-node  and included	 a central limit theorem for it. This time we will be  looking  at the degree of node $k$ for some   $k \in \mathbb{N}$.    We  begin by finding the first two moments of $d(k)$. 

\begin{theorem}\label{thm:expvardeggen}
Let $\mc{G}_n$ be  a random permutation graph on $n$ vertices, and let $d(k)$ be the degree of node $k \in \mathbb{N}$. Then, we have 
\begin{itemize}
\item[(i)] $\mbb{E}[d(k)] = \frac{n-1}{2}$; 
\item[(ii)] $Var(d(k)) = \frac{2(n - 1) + (n - 2k +1)^2}{12}$. 
\end{itemize}
\end{theorem}
\pf  Letting $X_i$'s be i.i.d. uniform random variables over $(0,1)$,  $d(k) =_d \sum_{j =1 }^{k-1} \mathbf{1}(X_j > X_k) + \sum_{j = k + 1 }^{n} \mathbf{1}(X_j <  X_k)$. 

(i) We have $\mathbb{E}[d(k)] = \sum_{j =1 }^{k-1} \mathbb{P}(X_j > X_k) + \sum_{j = k + 1 }^{n} \mathbb{P}(X_j <  X_k) = \frac{n-1}{2}$, since $X_i$'s are continuous random variables. 

(ii) Recall the law of total variance which in our case reads as 
\begin{equation}\label{vardkcomp1}
Var(d(k)) =  \mbb{E}[Var(d(k) \mid X_k)] +  Var(\mbb{E} [d(k) \mid X_k]).
\end{equation}
Beginning with the first term on the right-hand side of \eqref{vardkcomp1}, observing $$Var(d(k) \mid X_k) = \sum_{j=1}^{k-1} (1 - X_k) X_k + \sum_{j = k +1}^n  X_k (1 - X_k),$$  we have 
$$  \mbb{E}[Var(d(k) \mid X_k)] = \mbb{E} \left[ \sum_{j=1}^{k-1} (1 - X_k) X_k + \sum_{j = k +1}^n  X_k (1 - X_k)\right] = \frac{n-1}{6},$$ where we used $\mbb{E}[X_k] = 1 /2$ and $\mbb{E}[X_k^2] = 1 /3$. 

For the second term on right-hand side of \eqref{vardkcomp1},  note that $$\mbb{E}[d(k)\mid X_k] = \sum_{j=1}^{k-1} (1 - X_k) + \sum_{j=k+1}^n X_k = k - 1 + (n - 2k +1) X_k.$$ Therefore, $$ Var(\mbb{E} [d(k) \mid X_k]) = Var ( k - 1 + (n - 2k +1) X_k) = (n-2k +1)^2 \frac{1}{12}.$$ Using \eqref{vardkcomp1}, (ii) now follows. 
\bbox 

\begin{rmk}
(i) $\mathbb{E}[d(k)]$ is constant independent of $n$ due to symmetry.

(ii) Variance of $d(k)$ is minimized at $k = n /2$ and is maximized at   $k = 1$ and $k =n$. 
\end{rmk}

When we take $k = n /2$ in Theorem \ref{thm:expvardeggen}, we obtain $\mbb{E}[d(n /2)] = \frac{n-1}{2}$ and $Var(d(n/2)) = \frac{2n-1}{12}$. Given these, we may  prove an alternative central limit theorem for the mid-node with deterministic parameters in limit  in \eqref{eqn:cltwu}. We omit the details. 

How about other $k$'s? Next result is concerned with $k$ as a function of $n$ that is slowly growing. It applies in particular to the case where $k$ is fixed.

\begin{theorem}
Let $\mc{G}_n$ be  a random permutation graph on $n$ vertices and $h(n) = o (n)$ be an integer valued function. Then, $$\frac{d(h(n)) - n/2}{\sqrt{Var(d(h(n)))}} \rightarrow_d  \mc{Z}, \quad \text{as } n \rightarrow \infty,$$ with $Var(h(n)) = \frac{n^2 - n (1 + 4 h(n)) + 4 h(n) (h(n) +1)}{6} \sim \frac{n^2}{6}$. 
\end{theorem}

\pf 
  $X_i$'s are still  i.i.d. uniform random variables over $(0,1)$.  Let us use the convention $h = h(n)$. Observe
\begin{eqnarray}\label{eqn:hn1}
\nonumber d(h) &=_d& \sum_{j=1}^{h - 1} \mbf{1}(X_j > X_h) +  \sum_{j=h +1}^{n} \mbf{1}(X_j <  X_h) \\
\nonumber &=& \sum_{j=1}^{h - 1}  \mbf{1}(X_j > X_h)  + \sum_{j=h +1}^{n} (1 -  \mbf{1}(X_j > X_h) ) \\
\nonumber &=&  n -  h +   \sum_{j=1}^{h - 1}  \mbf{1}(X_j > X_h)  -  \sum_{j=h +1}^{n}   \mbf{1}(X_j > X_h) \\
&=& n  - \sum_{j = h+1}^{n-h}    \mbf{1}(X_j > X_h) - h +  \sum_{j=1}^{h - 1}  \mbf{1}(X_j > X_h)  - \sum_{j=n-h+1}^n  \mbf{1}(X_j > X_h). 
\end{eqnarray}
By assumption $h  = o(n)$, we have $h / n \rightarrow 0$. Note that the right-most two terms in \eqref{eqn:hn1} 
$$\frac{\sum_{j=1}^{h - 1}  \mbf{1}(X_j > X_h)  - \sum_{j=n-h}^n   \mbf{1}(X_j > X_h) }{n} \rightarrow 0,$$ with probability one as $n \rightarrow \infty$. This in particular implies the distributional convergence of $\frac{\sum_{j=1}^{h - 1}  \mbf{1}(X_j > X_h)  - \sum_{j=n-h}^n   \mbf{1}(X_j > X_h) }{n}$ to zero. 

Also, using symmetry techniques similar to ones in  previous section, the classical central limit theorem and the use of Slutsky's theorem imply that $$\frac{\frac{n}{2} - \sum_{j = h +1}^{n-h - 1} \mathbf{1}(X_j > X_k) }{n / \sqrt{6}}=_d \frac{\sum_{j = h +1}^{n-h - 1}  \mathbf{1}(X_j > X_k) - \frac{n}{2}}{n/ \sqrt{6}} \rightarrow_d \mc{Z}$$ as $n \rightarrow \infty$. With these      observations,  we conclude that 
\begin{eqnarray*} 
\nonumber \frac{d(h) - \frac{n}{2}}{n / \sqrt{6}} &=& 
\frac{ n  - \sum_{j = h+1}^{n-h}    \mbf{1}(X_j > X_h) - h +  \sum_{j=1}^{h - 1}  \mbf{1}(X_j > X_h)  - \sum_{j=n-h +1}^n  \mbf{1}(X_j > X_h) - \frac{n}{2}}{n / \sqrt{6}} \\
&& \longrightarrow_d \mc{Z} 
\end{eqnarray*}
as claimed.
\bbox

\bigskip

  The following corollary is now immediate. 

\begin{cor}
Let $\mc{G}_n$ be  a random permutation graph on $n$ vertices and $k \in \mbb{N}$ be fixed. Then, $$\frac{d(k) - n/2}{\sqrt{\frac{n^2 - 4nk + 2n + 4k^2 -2}{12}}} \rightarrow_d  \mc{Z}, \quad \text{as } n \rightarrow \infty. $$
\end{cor}
 
 \begin{remark}
One may wonder how the extremal degrees behave in a random permutation graph $\mc{G}_n$. Writing $\delta(\mc{G}_n)$ and $\Delta (\mc{G}_n)$  for the minimal and maximal degrees in $\mc{G}_n$, it was shown in  \cite{BM:2017} that  $\frac{\delta (\mc{G}_n)}{\sqrt{n}} \rightarrow_d \Gamma,$ as $n \rightarrow \infty$, where $\Gamma$  has a   Rayleigh distribution with parameter $1 / \sqrt{2}$, i.e., $\mbb{P}(\Gamma > \gamma) = e^{- \gamma^2}$ for all $\gamma > 0$. 
 The following  symmetry argument provides a similar result for  $\Delta(\mc{G}_n)$. 
Letting $X_1,X_2,\ldots$ be  a sequence of i.i.d. random variables uniformly distributed over $(0,1)$, we have  
\begin{eqnarray*}
\Delta(\mc{G}_n) &=_d& \max_{k \in [n]} \lp \sum_{j=1}^{k-1} \mathbf{1} (X_j > X_k)  +  \sum_{j=k+ 1}^{n} \mathbf{1} (X_j<  X_k)\rp \\ &=&  \max_{k \in [n]} \lp \sum_{j=1}^{k-1}( 1  -  \mathbf{1} (X_j <  X_k) )    +  \sum_{j=k+ 1}^{n} ( 1  -  \mathbf{1} (X_j >   X_k)) \rp  \\
&=&  \max_{k \in [n]}   \lp  n -  1 - \sum_{j=1}^{k-1} \mathbf{1} (X_j <  X_k)  - \sum_{j=k+ 1}^{n} \mathbf{1} (X_j>  X_k)\rp  \\
&=_d& n - 1 -   \min_{k \in [n]}   \lp \sum_{j=1}^{k-1} \mathbf{1} (X_j > X_k)  +  \sum_{j=k+ 1}^{n} \mathbf{1} (X_j<  X_k)\rp  \\
&=& n- 1  - \delta(\mc{G}_n).
\end{eqnarray*}
In particular, $\frac{n  - 1 - \Delta (\mc{G}_n)}{\sqrt{n}} \rightarrow_d \Gamma,$  as $n \rightarrow \infty$, where $\Gamma$ has  a   Rayleigh distribution with parameter $1 / \sqrt{2}$.
 	
\end{remark}
\section{Number of connected components}\label{sec:conncom}

In previous parts, we worked on the degree of a  specific node. Now  we will change the direction and look at the number of nodes with a given degree, for example, the number of isolated vertices, leaves, etc. In order to do so,  we will first have some discussion on connected components and probability of being connected.  The results of this section are somewhat available in the literature, we still include a detailed treatment for the sake of completeness. 
 
Before moving on to main discussion,  let us give some pointers to the literature.   First,  a permutation $\pi_n\in S_n$ is said to be \emph{decomposable} if there exists a positive integer $k<n$ such that  $\{\pi(1),\pi(2),\ldots,\pi(k)\}=\{1,2,\ldots,k\}$; otherwise, indecomposable.
Equivalently, a permutation $\pi_n\in S_n$ is indecomposable if $n$ is the least positive integer $k\in [n]$ for which $\{\pi(1), \pi(2), \ldots, \pi(k)\}=\{1,2, \ldots, k\} .$ We note that the corresponding graph $\mc{G}_n$ is connected if and only if $\pi_n$ is indecomposable.
Indecomposable permutations were probably considered for the first time by Lentin \cite{lentin} in 1970. Then Comtet  \cite{comtet} showed, letting $a_n$ be the number of indecomposable permutations of $[n],$ that 
\begin{equation}\label{eqn:gfindec}
\sum_{n \geq 1} a_n x^{n}=1-\frac{1}{\sum_{n \geq 0} n ! x^{n}}
\end{equation}
  Indecomposable permutations  also found their way into many textbooks on enumeration and combinatorics, see for example, \cite{fs} and  \cite{gouldjack}. 
 Lastly, let us note that the generating function  in \eqref{eqn:gfindec}  appear in a number of distinct areas. A look at  the  content of the entry A003319 in the On-line Dictionary Encyclopedia of Integer Sequences will make this clear (\url{https://oeis.org/A003319}).

Let us get back to our work  on the probability generating function for the probabilities of being connected for random permutation graphs. 
Throughout, let  $f_i$
 be the probability that a random permutation graph with $i$ vertices is connected. 
 
\begin{proposition} Let $\mc{G}_n$ be a random permutation graph on $n$ vertices. Then
\begin{equation}\label{eqn:exactk}
\Pro(\mc{G}_n \text{ has exactly } k \text{ components})=\sum_{b_1+b_2+\cdots+b_k=n, b_i \geq 1, i \in [k]}\prod_{j=1}^{k}f_{b_j}.
\end{equation}
\end{proposition}

\textbf{Proof.}
 Observe that $\mc{G}_n$ has exactly $k$ connected components if and only if there are $k$ vertices  $\{a_1,a_2,\ldots,a_k\}$ with $1 \leq a_1<a_2< \cdots <a_k = n$ such that
\begin{itemize}
\item $\{\pi(a_{j-1}+1),\pi(a_{j-1}+2),\ldots,\pi(a_j)\}=\{a_{j-1}+1,a_{j-1}+2,\ldots,a_{j}\}$, and 
\item $\{a_{j-1}+1,a_{j-1}+2,\ldots,a_{j}\}$ forms a connected subgraph,
\end{itemize} 
for all $j\in\{1,2,\ldots,k\}$. Here, we use the   convention $a_0=0$.  

Let us write $b_1=a_1$, $b_2=a_2-a_1$, $b_3=a_3-a_2$, \ldots, $b_k=a_k-a_{k-1}$. Note that $b_1+b_2+\cdots+b_k=n$ and $b_j\geq 1$ for each $j$. Since $\{a_{j-1}+1,a_{j-1}+2,\ldots,a_{j}\}$ forms a connected subgraph with probability $f_{b_j}$ for all $j$, the result follows.
\hfill $\square$

Writing $F(x)=\displaystyle{\sum_{n=1}^{\infty}f_nx^n}$, observe that  $\Pro(\mc{G}_n \text{ has exactly } k \text{ components})$ becomes the coefficient of $x^n$ in the  series $(F(x))^k$. We now focus on understanding the sequence $f_n$ a bit more starting with the following recursion.

\begin{proposition}  For $n\in\mathbb{N}$, we have 
$$\sum_{t=1}^{n}\dfrac{f_t}{\binom{n}{t}}=1. $$
\end{proposition}

\textbf{Proof.}
\noindent Let $q_k$ be the number of possible connected   permutation graphs on $k$ vertices, i.e. $q_k=f_k k!$. We will count the number of   permutation graphs on $n$ vertices in terms of the size of the connected component containing the vertex $1$. Denote   the size of the connected component containing the vertex $1$   by $t$. Then $\{\pi(1),\pi(2),\ldots,\pi(t)\}=\{1,2,\ldots,t\}$ with the condition that $\{1,2,\ldots,t\}$ forms a connected subgraph. Thus, there are   $q_t(n-t)!$   permutation graphs whose  size of the connected component containing $1$ is exactly $t$. 
Therefore,  $$n!=\sum_{t=1}^{n}q_t(n-t)!=\sum_{t=1}^{n}f_t t!(n-t)!=n!\sum_{t=1}^{n}\dfrac{f_t}{\binom{n}{t}},$$ from which our assertion follows. 
\hfill $\square$

As noted above, after proving this recursion in context of random permutation graphs, we realized that it also appears in various other work. See references provided above. 

 Having obtained a recursion for $f_n$, we  now obtain explicit formulas for these probabilities. 
In order to so we  find  the generating function for the  sequence $n! f_n $ and then read the coefficients.  

\begin{proposition}
The generating function for the sequence $n! f_n$ is given by 
$$Q(x)  = \sum_{m=1}^{\infty}(-1)^{m+1}\phi(x)^m,$$
 where $\phi(x)=  \ds\sum_{n=1}^{\infty} n! x^{n}$. \end{proposition}

\textbf{Proof.}
Recall that $\sum_{k=1}^{n}\dfrac{f_k}{\binom{n}{k}}=1$. This says $$f_n =1 - \sum_{k=1}^{n-1} \frac{f_k}{\binom{n}{k}}.$$ Multiplying both sides with $n!$ we obtain 
$$n!f_n =n! - \sum_{k=1}^{n-1} k! (n-k)! f_k.$$ 
Now multiplying both sides with $x^n$
 and summing over $n \geq 1$, we get 
\begin{eqnarray*}
 \sum_{n =1}^{\infty} n! f_n x^n &=& \sum_{n =1}^{\infty} n! x^n -  \sum_{n =1}^{\infty} \sum_{k=1}^{n-1} k!(n-k)!f_k x^n\\
   &=&\sum_{n =1}^{\infty} n! x^n  - \sum_{k=1}^{\infty} \sum_{n=k+1}^{\infty}  k!(n-k)!  f_k x^k x^{n-k} \\
  &=&\sum_{n =1}^{\infty} n! x^n  - \sum_{k=1}^{\infty} x^k k! f_k \sum_{n=k+1}^{\infty} (n-k)! x^{n-k} \\
 &=& \sum_{n =1}^{\infty} n! x^n  -\sum_{k=1}^{\infty} x^k k! f_k \sum_{n=1}^{\infty} n! x^{n} .
\end{eqnarray*} 
 Letting $\phi(x) = \ds\sum_{n=1}^{\infty} n! x^{n}$ and $Q(x)$ be the generating function for $n! f_n$'s, we obtain the relation $$Q(x) = \phi(x) (1 - Q(x)).$$ Solving it yields $$Q(x) = 1 - \frac{1}{1 + \phi (x)} = 1 - \frac{1}{1 - (-\phi (x))}  = 1-\sum_{m=0}^{\infty}(-1)^m\phi(x)^m = \sum_{m=1}^{\infty}(-1)^{m+1}\phi(x)^m$$
 \hfill $\square$

  Before stating the next  result, a piece of notation is in order. For  $n \in \mathbb{N}$, if $\lambda_1+\lambda_2+\ldots+\lambda_m=n$ where $\lambda_j$'s are positive integers, we write $\lambda=(\lambda_1,\lambda_2,\ldots,\lambda_m)\models n$ and $l(\lambda)=m$. 
 
\begin{theorem}
We have $$f_n = \dfrac{\ds\sum_{m=1}^{n}\sum_{\lambda\models n,l(\lambda)=m}(-1)^{m+1}\lambda_1!\lambda_2!\ldots\lambda_m!}{n!}$$
\end{theorem}

\textbf{Proof. }
Observe that $$f_n= \dfrac{1}{n!}\cdot [x^n]\sum_{m=1}^{\infty}(-1)^{m+1}\phi (x)^m =  \dfrac{1}{n!}\cdot [x^n]\sum_{m=1}^{\infty}(-1)^{m+1}\Big(\sum_{k=1}^{\infty}k!x^k\Big)^m$$
Notice that the coefficient of $x^n$ in $\ds\Big(\sum_{k=1}^{\infty}k!x^k\Big)^m$ equals to $\ds\sum_{\lambda\models n,l(\lambda)=m}\lambda_1!\lambda_2!\ldots\lambda_m!$ (which vanishes for $m>n$).  
\hfill $\square$

\begin{theorem}
 Let $\mc{G}_n$ be a random permutation graph on $n$ vertices.  We have $$\Pro(\mc{G}_n \text{ has exactly } k \text{ components})= \dfrac{1}{n(n-k)!}\sum_{w=k}^{n}\dfrac{\ds\sum_{\gamma\models n,l(\gamma)=w}(-1)^{w+k}\gamma_1!\gamma_2!\ldots\gamma_w!}{(w-1)(w-2)\cdots(w-k+1)}.$$
\end{theorem}
 
\textbf{ Proof.}
Recall the formula $$\Pro(\mc{G}_n \text{ has exactly } k \text{ components})=\ds\sum_{b_1+b_2+\ldots+b_k=n}\ds\prod_{j=1}^{k}f_{b_j}$$ from \eqref{eqn:exactk}. Let us write $q_{n,k}$ for the number of permutation graphs on $n$ vertices with exactly $k$ connected components. Now, for any $b=(b_1,b_2,\ldots,b_k)$ where $b\models n$, suppose $\lambda^{(j)}\models b_j$ with $l(\lambda^{(j)})=u_j$. Then, by putting $f_{b_j}$'s into this equation, we get $$q_{n,k}=\sum_{b\models n,l(b)=k}\prod_{j=1}^{k}\Big(\ds\sum_{u_j=1}^{n}\sum_{\lambda^{(j)}\models b_j,l(\lambda^{(j)})=u_j}(-1)^{u_j+1}\lambda^{(j)}_1!\lambda^{(j)}_2!\ldots\lambda^{(j)}_{u_j}!\Big)$$


Then, by rearranging terms and considering the double counting, we get
\begin{eqnarray*}
q_{n,k}&=&\ds \binom{n-1}{k-1}\sum_{w=k}^{n}\dfrac{\ds\sum_{\gamma\models n,l(\gamma)=w}(-1)^{w+k}\gamma_1!\gamma_2!\ldots\gamma_w!}{\binom{w-1}{k-1}}
\end{eqnarray*}
As a result, we find 
\begin{eqnarray*}
\Pro(\mathcal{G}_n\text{ has exactly }k\text{ components})&=&\ds\dfrac{\binom{n-1}{k-1}}{n!}\sum_{w=k}^{n}\dfrac{\ds\sum_{\gamma\models n,l(\gamma)=w}(-1)^{w+k}\gamma_1!\gamma_2!\ldots\gamma_w!}{\binom{w-1}{k-1}}\\
&=& \ds\dfrac{1}{n(n-k)!}\sum_{w=k}^{n}\dfrac{\ds\sum_{\gamma\models n,l(\gamma)=w}(-1)^{w+k}\gamma_1!\gamma_2!\ldots\gamma_w!}{(w-1)(w-2)\cdots(w-k+1)}
\end{eqnarray*}
\hfill $\square$


We conclude this section with some discussions on  $f_n$. We will be providing easy lower and upper bounds for it which will make a clear  connection to the number of isolated vertices. We start with the upper bound.

\begin{proposition}
For any $n \geq 1$, we have $$f_n \leq 1 - \frac{2}{n} + \frac{1}{n(n-1)}.$$
\end{proposition}

\textbf{Proof. }
It is easy to see that if $\pi_n(1) = 1 $ or $\pi_n(n) = n$, then the graph is disconnected. So 
\begin{eqnarray*}
 \mathbb{P}(\mathcal{G}_n \text{ is connected }) =  1 - \mathbb{P}(\mathcal{G}_n \text{ is disconnected}) &\leq& 1 - \mathbb{P}(\pi_n(1) = 1  \text{ or } \pi_n(n) = n) \\ &=& 1 - \frac{2}{n} + \frac{1}{n(n-1)}.
\end{eqnarray*}
\hfill $\square$

Next proposition provides an easy lower bound.

\begin{proposition}
For any $n \geq 1$, we have $$f_n \geq 1 - \frac{2}{n} + \Theta \left(\frac{1}{n^2}\right).$$
\end{proposition}

\textbf{Proof. }
$f_n =  \mathbb{P}(\mathcal{G}_n \text{ is connected}) = 1 -  \mathbb{P}(\mathcal{G}_n \text{ is disconnected})$. So let's find an upper bound for 
$\mathbb{P}(\mathcal{G}_n \text{ is disconnected})$. For $j=1,\ldots,n-1$, let $A_j $ be the event that $\{\pi(1), \ldots, \pi(j)\} = \{1,\ldots,j\}$. Then $\mathbb{P}(\mathcal{G}_n \text{ is disconnected}) \leq \mathbb{P} \left( \bigcup_{j=1}^{n-1} A_j\right)$. Note that $\mathbb{P}(A_1) = \mathbb{P}(A_{n-1}) = 1/n$. Also $\mathbb{P}(A_j) = \frac{j! (n-j)!}{n!}$. 
So 
\begin{eqnarray*}
\mathbb{P}(\mathcal{G}_n \text{ is disconnected}) &\leq& \mathbb{P} \left( \bigcup_{j=1}^{n-1} A_j\right) \\
&=& \mathbb{P}(A_1) + \mathbb{P}(A_{n-1}) + \mathbb{P}(A_2) + \mathbb{P}(A_{n-2}) + \sum_{j=3}^{n-3} \mathbb{P}(A_j) \\
&=& \frac{1}{n} + \frac{1}{n} + \frac{2}{n(n-1)} + \frac{2}{n(n-1)}  +  \sum_{j=3}^{n-3}  \frac{j! (n-j)!}{n!} \\
&=& \frac{2}{n} +  \frac{4}{n(n-1)} + \Theta \left(\frac{1}{n^2} \right) \\
&=&  \frac{2}{n} +  \Theta \left(\frac{1}{n^2} \right). 
\end{eqnarray*}
\hfill $\square$

Combining  the two observations we conclude that $$f_n = 1 - \frac{2}{n} + \Theta \left( \frac{1}{n^2} \right).$$
 Let us note that this was previously observed by Comtet    in a work related to indecomposable permutations. Indeed, he proves that a random permutation is indecomposable with probability $1 - 2 /n + E(n)$ where $E(n)$ is given by $$E(n)= -\frac{1}{(n)_{2}}-\frac{4}{(n)_{3}}-\frac{19}{(n)_{4}}-\frac{110}{(n)_{5}} 
-\frac{745}{(n)_{6}}-\frac{5752}{(n)_{7}}-\frac{49775}{(n)_{8}}-\frac{476994}{(n)_{9}}-\frac{5016069}{(n)_{10}}+O\Big(\frac{1}{n^{11}}\Big),$$
with  $(n)_{k}=n(n-1) \cdots(n-k+1)$, the falling factorial. However,  we decided to include the preceeding discussion due to its simple structure and due to the fact that it makes an  explicit relation between connectivity and isolated vertices   - below we will see that the expected number of isolated vertices in a random permutation graph behaves like $2/n$ as $n \rightarrow \infty$.

\section{Number of isolated vertices}\label{sec:numiso}

In this section we will  be in interested in the probability  $\Pro(\mc{G}_n \text{ has exactly } k \text{ isolated vertices})$.   Let us first look at the  permutation interpretation of  isolated vertices. 
  For a permutation  $\pi_n\in S_n$, $i$ is called a \emph{strong fixed point} of $\pi_n$ if  $\pi(j)<\pi(i) \text { for } 1 \leq j<i \text { and } \pi(j)>\pi(i) \text { for } i<j \leq n.$ 
Letting $I_n$ be the number of isolated vertices in a random permutation graph $\mc{G}_n$ with $n$ vertices, and  $\pi_n$  be the corresponding random permutation representation of $\mc{G}_n$, it is easily seen  that the isolated vertices in $\mc{G}_n$ are exactly the strong fixed points of $\pi_n$.

 Let below $i_n$   be  the probability that $\mc{G}_n$ has no isolated vertices. By convention we set   $i_0=1$ and $i_1=0$. Note that $i_n$ can also be considered as  the probability that each connected component of $\mc{G}_n$ has size at least two. Our main result in this section is Theorem \ref{thm:iso} which provides an explicit formula for $i_n$. To reach there we need some preparation.

Let $I(x)=\displaystyle{\sum_{t=0}^{\infty}i_tx^t}$ and $I_n$ be the number of permutation graphs on $n$ vertices having no isolated vertices, i.e., $I_n=n!i_n$. By conditioning  on the first connected component of $\mathcal{G}_n$, we get the recursion $$I_n=\ds\sum_{k=2}^{n}k!f_kI_{n-k}.$$ This yields:
\begin{proposition} 
We have 
$$i_n=\ds\sum_{k=2}^{n}\dfrac{f_ki_{n-k}}{\binom{n}{k}}. $$
\end{proposition}

Next result expresses  $\Pro(\mc{G}_n \text{ has exactly } k \text{ isolated vertices})$ in terms of  $i_n$'s.  

\begin{proposition}
We have 
\begin{equation}\label{kisoinnverfor}
\Pro(\mc{G}_n \text{ has exactly } k \text{ isolated vertices})=\sum_{b_0+b_1+\ldots+b_k=n-k,   b_i \geq 0, i \in [k]}\prod_{t=0}^{k}i_{b_t}.
\end{equation}
\end{proposition}

\textbf{Proof.}
Observe that the vertex $j$ is isolated in $\mc{G}_n$ if and only if the following hold:
\begin{itemize}
\item $\{\pi(1),\pi(2),\ldots,\pi(j-1)\}=\{1,2,\ldots,j-1\}$
\item $\{\pi(j+1),\pi(j+2),\ldots,\pi(n)\}=\{j+1,j+2,\ldots,n\}$
\item $\pi(j)=j$
\end{itemize}
Assume $\mc{G}_n$ has exactly $k$ isolated vertices, and let $S=\{j_1,j_2,\ldots,j_k\}$ be the set of isolated vertices. Similar to the calculations for  the number of connected components, let us write $b_0=j_1-1$, $b_1=j_2-j_1-1$, \ldots, $b_{k-1}=j_k-j_{k-1}-1$, $b_k=n-j_k$. Note that $b_0+b_1+\ldots+b_k=n-k$ and $b_t\geq 0$.  Given expression follows.
\hfill $\square$

\noindent Again, by using generating functions, $\Pro(\mc{G}_n \text{ has exactly } k \text{ isolated vertices})$ becomes the coefficient of $x^{n-k}$ in the power series $I(x)^{k+1}$.

Next result finds the generating function for the sequence $n! i_n$.

\begin{proposition}
The generating function of the sequence $n! i_n$ is 
$$J(x)=\sum_{n=0}^{\infty}n!x^n\cdot\sum_{m=0}^{\infty}(-1)^{m}\Big(\sum_{n=1}^{\infty}(n-1)!x^{n}\Big)^m.$$
\end{proposition}

\textbf{Proof.}
We have 
$$n!i_n=\sum_{k=2}^{n}k!f_ki_{n-k}(n-k)!$$
Now multiplying both sides by $x^n$
and summing over $n \geq 2$, we get 
\begin{eqnarray*}
\sum_{n=2}^{\infty}n!i_nx^n=\sum_{n=2}^{\infty}\sum_{k=2}^{n}k!f_ki_{n-k}(n-k)!x^n 
&=&\sum_{n=2}^{\infty}\sum_{k=2}^{n}k!f_kx^ki_{n-k}(n-k)!x^{n-k}\\
&=&\sum_{k=2}^{\infty}k!f_kx^k\sum_{n=k}^{\infty}i_{n-k}(n-k)!x^{n-k}\\
&=&\sum_{k=2}^{\infty}k!f_kx^k\sum_{n=0}^{\infty}i_nn!x^n
\end{eqnarray*}
Letting $J(x)=\ds\sum_{n=0}^{\infty}n!i_nx^n$, we have $J(x)-i_0-i_1x=(Q(x)-f_1x)J(x)$, which leads to $$J(x)-1=J(x)(Q(x)-x)\text{ and so }J(x)=\dfrac{1}{1+x-Q(x)}$$ 
Since $Q(x)=\ds\dfrac{\phi(x)}{1+\phi(x)}$, we can write $$J(x)=\dfrac{1+\phi(x)}{1+x+x\phi(x)}=\dfrac{\ds\sum_{n=0}^{\infty}n!x^n}{\ds1+\sum_{n=1}^{\infty}(n-1)!x^{n}}=\sum_{n=0}^{\infty}n!x^n\cdot\sum_{m=0}^{\infty}(-1)^{m}\Big(\sum_{n=1}^{\infty}(n-1)!x^{n}\Big)^m$$
\hfill $\square$

So, we have $i_n=\dfrac{1}{n!}[x^n]\ds\sum_{n=0}^{\infty}n!x^n\cdot\ds\sum_{m=0}^{\infty}(-1)^{m}\Big(\sum_{n=1}^{\infty}(n-1)!x^{n}\Big)^m$.  With the preceding notation, we get a closed form expression for $i_n$ as: 
\begin{theorem}
We have 
\begin{equation}\label{form:isvert}
i_n = \dfrac{n!+\ds\sum_{k=0}^{n-1}k!\Big(\sum_{m=1}^{n-k}\sum_{\lambda\models n-k,l(\lambda)=m}(-1)^m(\lambda_1-1)!(\lambda_2-1)!\ldots(\lambda_m-1)!\Big)}{n!}.
\end{equation}
\end{theorem}

 
The formula in next theorem for the probability of having exactly $k$ isolated vertices now follows from \eqref{form:isvert} combined with  a rewriting of  \eqref{kisoinnverfor}.

\begin{theorem}\label{thm:iso} We have
$$\Pro(\mc{G}_n \text{ has exactly } k \text{ isolated vertices})= \ds\sum_{\gamma\models n+1,l(\gamma)=k+1}\ds\prod_{s=1}^{k+1}i_{\gamma_s-1}, $$ where $$i_n = \dfrac{n!+\ds\sum_{k=0}^{n-1}k!\Big(\sum_{m=1}^{n-k}\sum_{\lambda\models n-k,l(\lambda)=m}(-1)^m(\lambda_1-1)!(\lambda_2-1)!\ldots(\lambda_m-1)!\Big)}{n!}.$$
In particular, the  expected number of isolated vertices is given by $\sum_{k=0}^n  k \ds\sum_{\gamma\models n+1,l(\gamma)=k+1}\ds\prod_{s=1}^{k+1}i_{\gamma_s-1}.$
\end{theorem}


\hfill $\square$

We will be dealing with  the  expected number of isolated vertices in a much simpler way in next section.

\section{Number of vertices with degree $d$}\label{sec:ndegd}

In previous section, we obtained the probabilities  that a random permutation graph has exactly $k$ isolated vertices, leading to a formula for the expected number of isolated vertices. The formula we derived for the expectation is not promising to get a simple expression.  The purpose of this section is to discuss an alternative approach to understand the mean number of isolated vertices. Indeed we will consider the more general case  of the number of vertices with degree $d$ for a given $d \geq 1$.

\begin{theorem}
Denote the number of vertices of degree $d$ in   $\mathcal{G}_n$ by $N_d$. Then $$\mathbb{E}[N_d] =\displaystyle{\sum_{j=1}^{n} \sum_{k=1}^{n}}\dfrac{\displaystyle{\binom{n-k}{\frac{d+j-k}{2}}\binom{k-1}{\frac{d+k-j}{2}}}}{n\cdot\displaystyle{\binom{n-1}{k-1}}}.$$
\end{theorem}

\textbf{Proof. } 	
Let $\pi_n$ be the corresponding random permutation of $\mc{G}_n$.We have 
\begin{eqnarray*}
	\mathbb{E}[N_d] &=& \sum_{j=1}^{n} \mathbb{E} [\mathbf{1}(\text{vertex j has degree d})] = \sum_{j=1}^n \mbb{P}(\text{vertex $j$ has degree $d$})  \\ &=& \sum_{j=1}^n \sum_{k=1}^n \mbb{P}(\text{vertex $j$ has degree $d$} \mid \pi_n(j)=k) \mathbb{P}(\pi_n(j)=k) \\
\end{eqnarray*}
Now, there are exactly $(n-1)!$ permutations in $S_n$ with $\pi_n(j)=k$ and $\mathbb{P}(\pi_n(j)=k)=\dfrac{1}{n}$ for fixed $j,k$. On the other hand, assume the vertex $j$ in $\mc{G}_n$ has degree $d$ where $\pi_n(j)=k$. Note that the vertex $i$ is adjacent to the vertex $j$ if and only if $i<j$ with $\pi_n(i)>k$ or $i>j$ with $\pi_n(i)<k$. If we denote the number of neighbors of $j$ in $\{1,2,\ldots,j-1\}$ by $s$, then there are exactly $s$ many elements in the intersection of $\{k+1,k+2,\ldots,n\}$ with $\{\pi_n(1),\pi_n(2),\ldots,\pi_n(j-1)\}$, which implies there are exactly $(k-1)-(d-s)$ many elements in the intersection of $\{1,2,\ldots,k-1\}$ with $\{\pi_n(1),\pi_n(2),\ldots,\pi_n(j-1)\}$. As a result, we get $s+(k-1)-(d-s)=j-1$ and so $s=\dfrac{d+j-k}{2}$.\\ 

Hence, given that $\pi_n(j)=k$, $\pi_n$ corresponds a permutation graph for which the vertex $j$ has degree $d$ if and only if the intersection of $\{k+1,k+2,\ldots,n\}$ with $\{\pi_n(1),\pi_n(2),\ldots,\pi_n(j-1)\}$ contains exactly $\dfrac{d+j-k}{2}$ many elements whereas the intersection of $\{1,2,\ldots,k-1\}$ with $\{\pi_n(1),\pi_n(2),\ldots,\pi_n(j-1)\}$ has cardinality $\dfrac{d+k-j}{2}$. Thus, we get $$\mbb{P}(\text{vertex $j$ has degree $d$} \mid \pi_n(j)=k)=\dfrac{\displaystyle{\binom{n-k}{\frac{d+j-k}{2}}\binom{k-1}{\frac{d+k-j}{2}}(k-1)!(n-k)!}}{(n-1)!}$$  \\

Therefore, we obtain $$\mathbb{E}[N_d]=\sum_{j=1}^{n} \sum_{k=1}^{n} \dfrac{\displaystyle{\binom{n-k}{\frac{d+j-k}{2}}\binom{k-1}{\frac{d+k-j}{2}}(k-1)!(n-k)!}}{n\cdot(n-1)!}=\displaystyle{\sum_{j=1}^{n} \sum_{k=1}^{n}}\dfrac{\displaystyle{\binom{n-k}{\frac{d+j-k}{2}}\binom{k-1}{\frac{d+k-j}{2}}}}{n\cdot\displaystyle{\binom{n-1}{k-1}}}.$$
\hfill $\square$

   \begin{corollary}
   The expected number of isolated vertices $I_n$ is given by $$\mathbb{E}[I_n]= \frac{1}{n} \sum_{k=1}^n \frac{1}{\binom{n-1}{k-1}}$$
   \end{corollary}
  
  In order to understand this last expression, we need a computational proposition. 

\begin{propn}
Let $n \in \mbb{N}$. Then $$2 + \frac{2}{n} \leq  \sum_{k=0}^{n} \frac{1}{\binom{n}{k}} \leq  2 + \frac{4}{n}.$$    In particular, $$\limm_{n \rightarrow \infty}  \sum_{k=0}^{n} \frac{1}{\binom{n}{k}} = 2.$$
\end{propn}

\pf
Note that  when $2 \leq k \leq n-2$, we have $\frac{1}{\binom{n}{k}} \leq \frac{1}{\binom{n}{2}} = \frac{2}{n (n-1)}$. But then $$1 + \frac{1}{n} + \frac{1}{n} + 1 \leq \sum_{k=0}^{n} \frac{1}{\binom{n}{k}}  \leq  1 + \frac{1}{n} + \frac{1}{n} + 1  + \sum_{k=2}^{n-2} \frac{2}{n(n-1)}.$$ So $$2 + \frac{2}{n} \leq  \sum_{k=0}^{n} \frac{1}{\binom{n}{k}}  \leq  2 + \frac{2}{n} + \frac{2 (n-3)}{n (n-1)} \leq  2 + \frac{4}{n}.$$ The limit result   follows immediately. 
\bbox

\begin{corollary}
Let $I_n$ be the number of isolated vertices   $\mc{G}_n$. Then $$\frac{n \mbb{E}[I_n]}{2} \rightarrow 1, \quad \text{as } n \rightarrow \infty.$$ 
\end{corollary}

Lastly we give the expectation of the number of leaves. 

   \begin{corollary}
   The expected number of leaves is given by $$\mathbb{E}[N_1] = \frac{2}{n} \sum_{j=1}^n \frac{j}{\binom{n-1}{j}}.$$ $\mathbb{E}[N_1]  \sim 2$ as $n \rightarrow \infty$.
   \end{corollary}

\textbf{Proof.} 
We have 
\begin{eqnarray*}
\mathbb{E}[N_1]  = \sum_{j=1}^n \sum_{|k - j|=1} \dfrac{\displaystyle{\binom{n-k}{\frac{1+j-k}{2}}\binom{k-1}{\frac{1+k-j}{2}}}}{n\cdot\displaystyle{\binom{n-1}{k-1}}} &=& \frac{1}{n} \sum_{j=1}^n \frac{j}{\binom{n-1}{j}} +    \frac{1}{n} \sum_{j=1}^n \frac{n-j + 1}{\binom{n-1}{j- 2}} \\
&=& \frac{2}{n} \sum_{j=1}^n \frac{j}{\binom{n-1}{j}}.
\end{eqnarray*}
For the asymptotics  note that 
$$ \frac{2}{n} \sum_{j=1}^n \frac{j}{\binom{n-1}{j}} = \frac{2}{n} \left(\frac{1}{n - 1} + \frac{4}{n(n-1)} + \cdots + \frac{n-2}{n-1} + n-1 \right) \sim 2$$
 as $n \rightarrow \infty$.
\hfill $\square$

It is likely that the number of leaves in a random permutation has a Poisson behavior, but we were not able to justify this yet.

\section{Number of $m$-cliques}\label{sec:clique}

Let $K_m$ be the number of $m$-cliques, $m \in \mathbb{N}$,  in a random permutation graph $\mc{G}_n$ with $n$ vertices. Let $\pi_n$ be the corresponding permutation representation in $S_n$.  Clearly, $K_1 = n$ and $K_2$ is the number of edges in $\mc{G}_n$ which of course is just the number of inversions in $\pi_n$, denoted by  $Inv(\pi_n)$.     The permutation statistic $Inv(\pi_n)$ is well-studied in the literature, and there are various different proofs that it satisfies a central limit theorem when $\pi_n$ is uniformly random. Let us just refer to  \cite{fulman} where not only $\frac{Inv(\pi_n) - \binom{n}{2}}{\sqrt{\frac{n(n-1) (2n +5 )}{72}}}  \rightarrow_d \mc{Z}$ as $n \rightarrow \infty$ is proven, but it is also shown that a convergence rate of order $1 \sqrt{n}$ can be obtained with respect to $d_K$.   In this section, we focus on $K_m$ for a given $m \geq 2$ and show that it satisfies a central limit theorem.   Our main tools will be the random permutation interpretation and the theory of $U$-statistics. 

\begin{theorem}\label{thm:Km}
Let  $K_m$ be the number of $m$-cliques in  $\mc{G}_n$. We then have $$\frac{K_m - \mathbb{E}[K_m] }{\sqrt{Var(K_m)}}
\longrightarrow_d \mathcal{Z}, \qquad \text{as} \; \; n\rightarrow \infty,$$ where $$\mbb{E}[K_m] = \binom{n}{m} \frac{1}{m!},$$ and where   $Var(K_m) = \mbb{E}[K_m^2] - \lp \mbb{E}[K_m] \rp^2$, with \begin{equation}\label{eqn:2ndmomentKm}
\mathbb{E}[K_m^2]= \sum_{t+s \leq m} [(2 m - t)!]^{-1} 4^{m-t} \binom{n}{2m-t}\binom{m-t-s-1/2}{m-t-s} \binom{s+(t+1)/2-1}{s} \binom{2m-t}{2m-2t-2s}.
\end{equation}
\end{theorem}
\begin{remark}
 It is further shown in \cite{IO:2018} that  $$Var(K_m) \sim_{n \rightarrow \infty} \frac{1}{2 ((2 m - 1)!)^2}  \lp \binom{4m - 2}{2m -1} - 2 \binom{2m - 1}{m}^2 \rp n^{2m - 1},$$ therefore an application of Slutsky' theorem yields a slightly more compact form 
$$\frac{K_m  - \binom{n}{m} \frac{1}{m!}}{ \frac{1}{2 ((2 m - 1)!)^2}  \lp \binom{4m - 2}{2m -1} - 2 \binom{2m - 1}{m}^2 \rp n^{2m - 1}} \rightarrow_{d} \mc{Z}, \quad n \rightarrow \infty.$$
\end{remark}

\begin{proof}
Let $\pi_n$ be the corresponding random permutation to $\mc{G}_n$. We first start by observing that a subset $\mc{S} = \{i_1, i_2,\ldots,i_m\}$, $i_1< i_2 < \cdots < i_m$ forms an $m$-clique in $\mc{G}_n$ if and only if $$\pi_{i_1} > \pi_{i_2} > \cdots > \pi_{i_m}.$$ Using this, we  then write $$K_m = \sum_{1 \leq i_1 < i_2 < \cdots < i_m \leq n} \mathbf{1}(\pi_{i_1} > \pi_{i_2} > \cdots > \pi_{i_m}).$$ In other words, $K_m$ is merely the number of decreasing subsequences of length $m$ in $\pi_n$. Noting that this equals in distribution to the number of increasing subsequences of length $m$ in $\pi_n$, and denoting the latter by $I_{n,m}$, \cite{IO:2018} shows that $\mbb{E}[I_{n,m}] = \binom{n}{m} \frac{1}{m!}$  and  $Var(I_{n,m}) = \mbb{E}[I_{n,m}^2] - \lp \mbb{E}[I_{n,m}] \rp^2$  where $\mbb{E}[I_{n,m}^2]$ is as given in \eqref{eqn:2ndmomentKm}.  
From these observations, the moments given in the statement of Theorem \ref{thm:Km} are clear.\\

Now we move on to proving the central limit theorem. Although a brief sketch  for the number of  increasing subsequences of a given length  in a random permutation  was given in \cite{IO:2018}, we include all the details here  for the sake of completeness\footnote{The reason for why it was only a sketch in \cite{IO:2018} is that a similar argument on random words was detailed there.}. First, the proof is based on $U$-statistics, and   a general reference for such statistics is \cite{KT:2008}.  Remember that   $$K_m = \sum_{1 \leq i_1 < i_2 < \cdots < i_m \leq n} \mathbf{1}(\pi_{i_1} > \pi_{i_2} > \cdots > \pi_{i_m})=_d \sum_{1 \leq i_1 < i_2 < \cdots < i_m \leq n} \mathbf{1}(\pi_{i_1} < \pi_{i_2} < \cdots < \pi_{i_m}).$$ Letting then  $X_1,\ldots,X_n$ be i.i.d.
random variables from some continuous distribution, a result often attributed to R\'enyi tells us that  
$$ \sum_{1 \leq i_1 < i_2 < \cdots < i_m \leq n} \mathbf{1}(\pi_{i_1} < \pi_{i_2} < \cdots < \pi_{i_m}) =_d  \sum_{1 \leq i_1 < i_2 < \cdots < i_m \leq n} \mathbf{1}(X_{i_1} < X_{i_2} < \cdots < X_{i_m}).$$
 Let now $Y_1,\ldots,Y_n$ be another i.i.d. sequence of random variables from some continuous distribution. Assume further that the families $\{X_1,\ldots,X_n\}$ and $\{Y_1,\ldots,Y_n\}$ are independent as well. Define the (random) permutation $\gamma \in S_n$ with the property  $$Y_{\gamma(1)}< Y_{\gamma(2)} < \cdots < Y_{\gamma(n)}.$$ Clearly, the i.i.d. assumption yields $(X_1,\ldots,X_n)  =_d (X_{\gamma(1)},\ldots,X_{\gamma(n)})$. But then, letting  $\mc{R} = \{(i_1,\ldots,i_m): \text{each } i_j \in [n], i_j\text{'s  are distinct from each other}\}$
\begin{eqnarray*} \\
  K_m &=_d&   \sum_{1 \leq i_1 < i_2 < \cdots < i_m \leq n} \mathbf{1}(X_{i_1} < X_{i_2} < \cdots < X_{i_m}) \\ &=_d&   \sum_{1 \leq i_1 < i_2 < \cdots < i_m \leq n}  \mathbf{1}(X_{\gamma(i_1)} < X_{\gamma(i_2)} < \cdots < X_{\gamma(i_m)}) \\
  &=&  \sum_{\mc{R}}  \mathbf{1}(X_{\gamma(i_1)} < X_{\gamma(i_2)} < \cdots < X_{\gamma(i_m)}, i_1 < i_2 < \cdots < i_m) \\
  &=&    \sum_{\mc{R}}  \mathbf{1}(X_{\gamma(i_1)} < X_{\gamma(i_2)} < \cdots < X_{\gamma(i_m)}, Y_{\gamma(i_1)} <Y_{\gamma(i_2)} < \cdots < Y_{\gamma(i_m)})  \\
  &=_d&  \sum_{\mc{R}}  \mathbf{1}(X_{i_1} < X_{i_2} < \cdots < X_{i_m}, Y_{i_1} <Y_{i_2} < \cdots < Y_{i_m}).
\end{eqnarray*}
Here, in the last step we used the fact that $\gamma$ is a bijection, being a member of the symmetric group. 

 Now we define the function $$h((x_{i_1},y_{i_1}),\ldots,(x_{i_m},y_{i_m})) = \sum f((x_{j_1},y_{j_1}),\ldots,(x_{j_m},y_{j_m})),$$ where the summation is over all permutations of $i_1,\ldots,i_m$ and where $$ f((x_{i_1},y_{i_1}),\ldots,(x_{i_m},y_{i_m})) = \frac{1}{\binom{n}{m}}  \mathbf{1}(x_{i_1} < x_{i_2} < \cdots < x_{i_m}, y_{i_1} <y_{i_2} < \cdots < y_{i_m}).$$ Note that $h$ is clearly symmetric.  Therefore, we are able to express $$K_m =_d \frac{1}{\binom{n}{m}} \sum_{1 \leq i_1 < i_2 < \cdots < i_m \leq n} h((X_{i_1}, Y_{i_1}),\ldots,(X_{i_m}, Y_{i_m}))$$ But then $K_m$ is a U-statistic because 
\begin{itemize}
\item[i.] $h$ is symmetric;
\item[ii.] $h \in L^2$ for each $n$, since it is just a finite sum of  indicators;
\item[iii.] $h$ is a function of independent random variables (or, vectors, to be more precise). 
\end{itemize} 
 Result follows. 
\hfill $\square$
\end{proof}

Let $C_m$ be the number  cycles of size at least  $m$ in a random permutation graph corresponding to the random permutation $\pi_n$. Three observations: 
\begin{enumerate}
\item $C_m$ is equal in distribution to the number of decreasing subsequences in a random permutation of length at least $m$.  
\item The distribution of the number of decreasing subsequences in $\pi_n$ is the same as the distribution of the number of increasing subsequences in $\pi_n$; 
\item  Letting $X_1,X_2,\ldots$ be i.i.d. random variables uniform over the interval $(0,1)$ $$C_m =_d \sum_{1 \leq i_1 < i_2 \leq \cdots < i_m \leq n} \mbf{1}(X_{i_1}< X_{i_2} < \cdots < X_{i_m}).$$ 
\end{enumerate}
But we are already familiar with the the last statistic on the right-hand side  from our discussion on $m$-cliques, and the exact  arguments there show that: 

\begin{theorem}
For a sequence of  random permutation graphs $\mc{G}_n$,  the number of cycles of length at least $m $ satisfies the central limit theorem:  $$\frac{C_m  - \binom{n}{m} \frac{1}{m!}}{ \frac{1}{2 ((2 m - 1)!)^2}  \lp \binom{4m - 2}{2m -1} - 2 \binom{2m - 1}{m}^2 \rp n^{2m - 1}} \rightarrow_{d} \mc{Z}, \quad n \rightarrow \infty.$$
\end{theorem}
 
Let us note that the literature on the number of  increasing subsequences of a random permutation is vast.  We will next mention a few more pointers how these interpret in terms  of cycles of   random permutation graphs. As a first result, if we are not only  interested in cycles of given sizes, but in all cycles of a random permutation graph, then the results of  Lifschitz and Pittel  \cite{LP:1981} are really useful. Letting $C_0 = 1$ by definition,  denoting the number of all cycles in $\mc{G}_n$ by $C_n^*$, and interpreting   results of  Lifschitz and Pittel   in our setting, we obtain 
\begin{equation*}\label{eqn:ECn}\mathbb{E}[C_n^*] =
\sum_{m=0}^n \frac{1}{m!} \binom{n}{m}
\end{equation*} and
\begin{equation*}\label{eqn:VarCn}\mathbb{E}[(C_n^*)^2] = \sum_{m+\ell \leq n} 4^{\ell} ((m+\ell)!)^{-1} \binom{n}{m+\ell} \binom{(m+1)/2+\ell-1}{\ell}.
\end{equation*}
Moreover, in the same paper, they prove  certain asymptotic relations which in our case may be read as  
$$\mathbb{E}[C_n^*] \sim (2 \sqrt{\pi e} )^{-1} n^{-1/4}
\exp(2n^{1/2}),$$ and
$$\mathbb{E}[(C_n^*)^2] \sim cn^{-1/4} \exp\left(2 \sqrt{2 + \sqrt{5}}
n^{1/2}\right)$$  as  $n \rightarrow \infty$, where $c \approx 0.0106$. It is not hard to see that with these moment asymptotics, a classical central limit theorem does not hold for $C_n^*$ which is slightly disappointing, but of course maybe there is some other distributional convergence which we do not know yet.

Another question that could be understood here is  the length of the longest cycle in a given random permutation graph $\mc{G}_n$. Let us denote this statistic by  $L_n = L_n(\mc{G}_n)$ By our discussions above, if we let $X_1,X_2,\ldots,$ be an i.i.d. sequence of uniform random variables over $(0,1)$, it should be clear that  $L_n $  has the same distribution as  the largest $k$ so that there exists some $1 \leq i_1 < i_2 \leq \cdots < i_k \leq n$  that satisfies $X_{i_1}< X_{i_2} < \cdots < X_{i_k}$. But then our problem just reduces to the standard longest increasing subsequence problem which was  solved by Baik, Deift and Johansson \cite{BDJ:1999}. Their results imply that 

\begin{theorem} The length of the longest cycle $L_n$ of  $\mc{G}_n$ satisfies 
$\frac{L_n - 2 \sqrt{n}}{n^{1/6}}
\longrightarrow_d TW$ as  $n \rightarrow \infty$,
 where $TW$ is the Tracy-Widom distribution whose cumulative distribution function is given by  $F(t) =
\exp\left(-\int_t^{\infty} (x-t) u^2(x) dx \right)$ where $u(x)$ is the solution of the Painlev\'e II equation  $u_{xx}=2u^3+xu \qquad \text{with} \qquad
u(x) \sim -Ai(x)\quad \text{as} \quad x \rightarrow \infty, $ and  $Ai(x)$ is the Airy's function.
\end{theorem}

\bigskip

\noindent \textbf{Acknowledgements:}  The second author is supported by the Scientific and Research
Council of Turkey [TUBITAK-117C047].  We would like to thank T{\i}naz Ekim A\c{s}{\i}c{\i} for introducing us random permutation graphs. Also, parts of this paper were completed at the Nesin Mathematics Village, the authors would like to thank Nesin Mathematics Village for their kind hospitality.

 \end{document}